\documentclass[11pt]{article}
\author{E. Mac\'{\i}as-Virg\'os}
\title{Length of parallel curves}
\usepackage{times}
\usepackage{graphicx}
\usepackage{url}
\usepackage[usenames]{color}

\date{}

\newtheorem{prop}{Proposition}

\newtheorem{theorem}[prop]{Theorem}
\newtheorem{defi}{Definition}
\newtheorem{cor}[prop]{Corollary}

\newcommand{\normal}{\overrightarrow{\bf n}}
\newcommand{\tangent}{\overrightarrow{\bf t}}
\newcommand{\exterior}{\overrightarrow{\bf e}}

\begin{document}
\maketitle
 
\begin{abstract}
We prove that the length difference  between a closed periodic curve   and its parallel curve at a sufficiently small distance $\varepsilon$ equals $2\pi\varepsilon$ times the rotation index. As an application, the rotation index of a curve could be estimated by means of  Cauchy-Crofton's formula.
\end{abstract}

\paragraph*{INTRODUCTION.}


The aim of this   note is to prove the following result. Let $\alpha$ be a closed periodic regular curve, let $\beta$ be the  parallel curve at distance $\varepsilon\geq0$. Assume that $\varepsilon$ is small enough  to not exceed the radius of curvature  of $\alpha$ when  $\kappa>0$ ($\kappa$ is the signed curvature of $\alpha$). Let $\omega$ be the rotation index of $\alpha$.

 \begin{theorem}\label{MAIN}  The length difference  $L(\alpha)-L(\beta)$ equals $2\pi\varepsilon\omega$. 
 \end{theorem} 
 
 I think this result is new or at least is not well known in differential geometry. Although elementary, it seems interesting because it actually finds the exact difference and shows that a relatively sophisticated invariant like the rotation
index can be determined by a much simpler invariant, namely, the length of a curve. Of
course, the converse is also a useful observation: computing the length by the rotation index
since this is a regular homotopy invariant \cite[p.330]{MONTIELROS}. As a corollary, the difference of the length of a
curve and its $\varepsilon$-parallel curve is a regular homotopy invariant. Possibly this could also be
directly proved by a variational argument..

\paragraph*{BASIC DEFINITIONS AND NOTATIONS.} 
Let $\alpha(t)$ be a differentiable plane curve, defined on the interval  $[a,b]$. The {\em length} of the curve is given by
\begin{equation}\label{LONGITUD}L(\alpha)= \int_a^b{\vert \alpha^\prime(t)\vert\,dt}.
\end{equation}
Suppose that the curve is regular, which means that the speed  $\vert\alpha^\prime \vert$ never vanishes.
Then the {\it arc-length parameter}  $s(t)$,  defined by $ds=\vert \alpha^\prime \vert dt$ and $s(a)=0$,
serves to reparametrize the curve with unit speed.

The (signed) {\em curvature} of $\alpha$ is the  function
\begin{equation}\label{CURVATURA}
\kappa = \det( \alpha^\prime ,\alpha^{\prime\prime} )/ \vert \alpha^\prime \vert^3.
\end{equation}
If the parameter is arc-length, the absolute value of the curvature is $\vert \kappa \vert =\vert \ddot{\alpha}\vert$, the module of the second derivative. 

When $\kappa \neq 0$, 
the {\em unitary normal vector} $\normal=\ddot{\alpha}/\vert \ddot{\alpha}\vert$ is well defined. It is perpendicular to the  tangent direction  and it points inwards the curve. For an arbitrary parameter $t$, the vector $\alpha^{\prime\prime}$ is not collinear to $\ddot{\alpha}$,  but both  are on the same side of the tangent line.

\paragraph*{PARALLEL CURVES.} Let $\alpha(t)$ be an arbitrary regular parametrization.  
 At each point $\alpha(t)$ we  take a unitary vector $\overrightarrow{{\bf e}(t)}$ orthogonal to $\alpha^\prime(t)$ and such that $\det(\alpha^\prime,\exterior)>0$. In other words, $\exterior$ is obtained by rotating in the  counter-clockwise sense the unitary tangent vector $\tangent=\dot{\alpha}=\alpha^\prime/\vert\alpha^\prime\vert$ (in Alfred Gray's book \cite{GRAY}, $\exterior$ is denoted by $J\alpha^\prime$).
 
 Then $\exterior  = + \normal $ when  $\kappa>0$ (the curve turns left) and  $\exterior  = - \normal$ when $\kappa<0$ (the curve turns right).
 
 \begin{defi}
 We define the (left) {\em parallel} curve to $\alpha$ at distance $\varepsilon \geq 0$ as the curve $\beta=\alpha+\varepsilon \exterior$.
 \end{defi}
 
 {\it Remark:\/} It is unnecessary to consider the case $\varepsilon\leq 0$, as we can always reparametrize the curve
 $\alpha$ in the opposite direction.\\
 \begin{figure}[htbp]
\begin{center}
 \includegraphics[width=4cm, height=4cm]{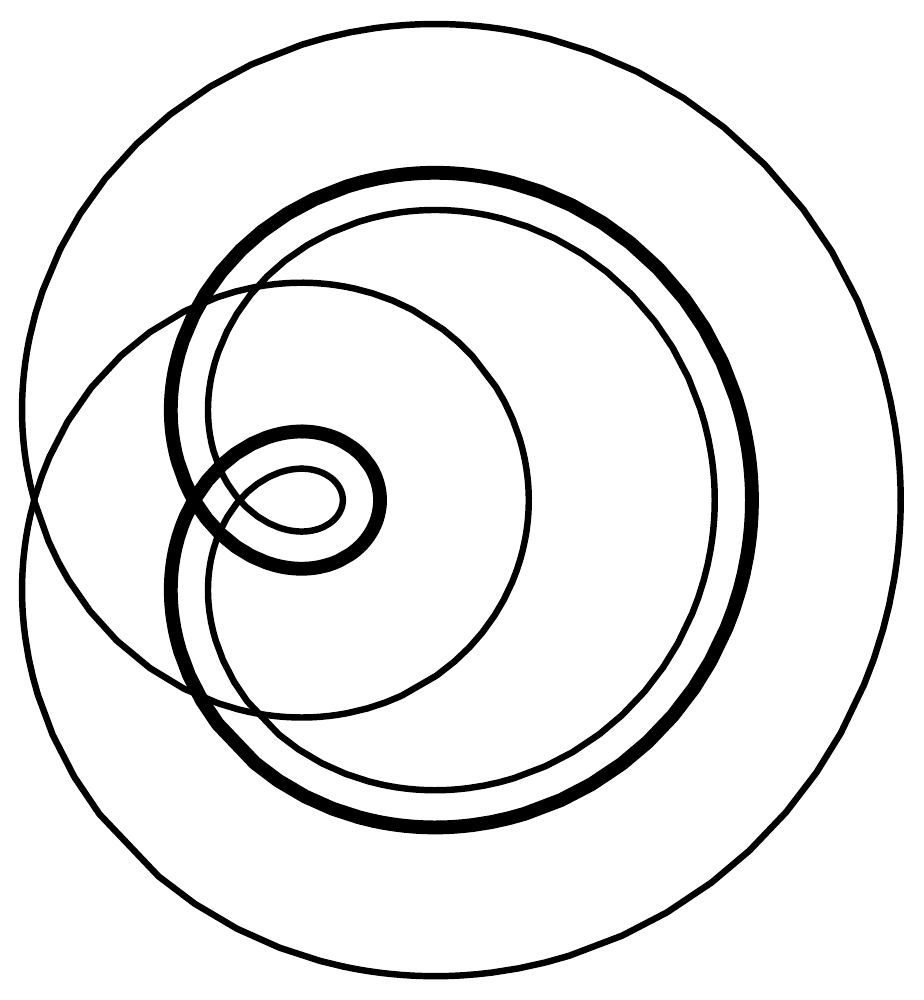}
 \caption{A curve (Pascal Snail) and two parallels}
\label{FIGURE1}
\end{center}
\end{figure}
We now   discuss the regularity of   $\beta $. For that we have to  take into account the {\em radius of curvature} $1/\vert \kappa \vert$ and the {\em evolute} of $\alpha$, which is the geometric locus of the {\em centers of curvature} $\alpha+(1/\vert\kappa\vert)\normal = \alpha +(1/\kappa)\exterior$. 

By differentiating with respect to the arc-length parameter $s$ of $\alpha$, we obtain
$\dot{\exterior }= -\kappa \tangent$, which is just a reformulation of the usual 
{\it Fr\'enet formula} $\dot{\normal}=-\vert\kappa\vert\tangent$  \cite{MONTIELROS,DOCARMO,GRAY}.
Hence $d\beta/ds=(1 -\varepsilon  \kappa)\dot{\alpha}$ and
\begin{equation}\label{SING}
 \vert d\beta/ds \vert = \vert 1 -\varepsilon  \kappa\vert.
\end{equation}

\begin{figure}[htbp]
\begin{center}
 \includegraphics[width=4cm, height=4cm]{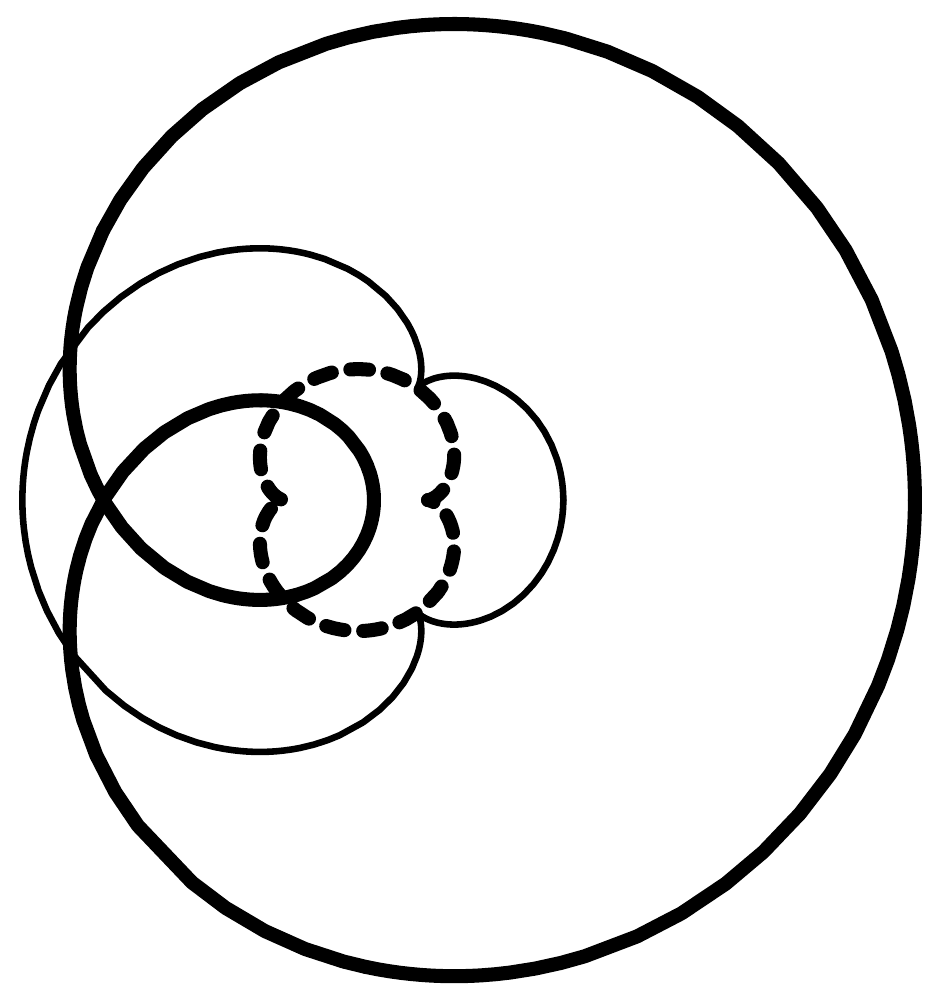}
 \caption{The same curve, its evolute (dashed) and one parallel with two singularities}
\label{FIGURE2}
\end{center}
\end{figure}

It follows that the parallel $\beta$ has a singularity each time $\varepsilon$ equals $1/\kappa$. This can only occur (as we are taking $\varepsilon\geq 0$) when $\kappa>0$ and $\varepsilon$ equals the radius of curvature, i.e. the parallel $\beta$ touches the evolute of $\alpha$ at corresponding points (see figure \ref{FIGURE2}).\\

{\it Remark:} The evolute itself has singularities at the places where the curvature attains a critical value; this is a consequence of the fact that the tangent vector to the evolute points in the {\em normal} direction to $\alpha$. \\
By  applying definition (\ref{LONGITUD}) to formula (\ref{SING}) we obtain the length of $\beta$.

\begin{theorem}\label{LONG} The length of the left parallel curve $\beta$ at distance $\varepsilon\geq 0$ to $\alpha$ is given by
$$L(\beta) = \int_0^{L(\alpha)}{\vert 1-\varepsilon \kappa(s)\vert \,ds}.$$
\end{theorem}

In Corollary \ref{EXTERIOR} we shall emphasize two particular cases of Theorem \ref{LONG}.

\begin{defi} The   {\em total curvature} of the curve $\alpha$ is the number
$$K=\int_0^{L(\alpha)}{\kappa(s)\,ds} = \int_a^b{\kappa(t)\vert \alpha^\prime(t)\vert\,dt}.$$
\end{defi}

\begin{cor}\label{EXTERIOR}
\begin{enumerate}
\item
 If $\kappa \leq 1/\varepsilon$ then 
$L(\beta)=L(\alpha)-\varepsilon  K;$
\item
 If $\kappa\geq0$ and $\varepsilon\geq  1/ \kappa  $ then 
$L(\beta)= \varepsilon  K -  L(\alpha)$.
\end{enumerate}
\end{cor}

{\it Example:} Let  $\alpha(t)=(R\cos t,R\sin t)$, $0\leq t \leq \pi$, be a half-circle with a big radius $R>0$. It has global curvature $K=\pi$. The parallel curve at distance $R+1$ is a small half-circle of radius $1$ wich goes backwards. Its length is $(R+1) \pi   - \pi R =\pi.$  \\


{\em Remark:\/}  From (\ref{CURVATURA}) it follows that the curvature of the parallel curve $\beta=\alpha+\varepsilon \exterior$  is given by
$\kappa_\beta={\kappa  / \vert 1-\varepsilon\kappa \vert}$, see   \cite[p. 117]{GRAY}. Then, when $\kappa < 1/\varepsilon$,   $\beta$   has the same evolute that $\alpha$.\\

\paragraph*{CLOSED CURVES} 
Let $\alpha(t)$  be  a regular curve defined in $[a,b]$. 
From now on we shall assume that our  curve is {\em closed}  and {\em periodic}, i.e.  it satisfies $\alpha(a)=\alpha(b)$ and
$\alpha^\prime(a)=\alpha^\prime(b).$\\

Let us recall the notion of {\it rotation index} (also called {\it turning number}). For simplicity, we parametrize $\alpha$ by the arc-length $s\in[0,L(\alpha)]$,
so the tangent vector $\tangent = \dot\alpha$ has module $1$. Write $\dot\alpha = (\cos \theta, \sin \theta)$. Then
$$\kappa = \det(\dot\alpha,\ddot\alpha)=d\theta/ds,$$
which proves that it is always possible to choose the angle $\theta$ in a  differentiable way (unique for any preassigned value of $\theta(0)$). Namely
\begin{equation}\label{ANGULO}
\theta(s)=\theta(0)+\int_0^s{\kappa}.
\end{equation}

Clearly $\theta$ does not depend on the parametrization of $\alpha$. Moreover, since the curve is  periodic,
the difference $\theta(b)-\theta(a)$ equals $2\pi \omega$,  for some  integer number $\omega$.

\begin{defi}[\mbox{\cite[p. 159]{GRAY}}]   The integer $\omega$ is called the {\em rotation index} of $\alpha$. It measures how many times the curve turns with respect to a fixed direction.
\end{defi}

{\it Example:\/} The rotation index of the Pascal Snail in Figure \ref{FIGURE1} is $\omega = \pm 2$  { depending on the sense of rotation}.\\

The following result is immediate from (\ref{ANGULO}).

\begin{prop}\label{TIERRA} The 
total curvature $K$ of a closed periodic curve with rotation index $\omega$ equals $2\pi \omega$.
\end{prop}

Finally,  if $\varepsilon$ is not too large, directly from Corollary \ref{EXTERIOR} we obtain Theorem \ref{MAIN}.

In addition, we have the following consequence of the Hopf theorem on turning tangents \cite[p. 333]{MONTIELROS}.

\begin{cor} For a simple closed curve in the plane (i.e., one without selfintersections), the
length of the $\varepsilon$-parallel curve minus the length of the original curve is always $\pm2\pi\varepsilon$, for $\varepsilon$ small enough.
\end{cor}
 


{\it Remark:\/} The following is a very well-known fact, which seems quite striking to non\--ma\-the\-ma\-ticians. Imagine that we surround the earth by the ecuator with a cable at ground level. If we next wanted the cable to stand a metre above ground level, how much extra cable would we need? The answer is: a little more than 6 metres. The reason is:  $2\pi(R+1)-2\pi R = 2\pi$. Of course this is a very particular case of our result.


\paragraph*{ESTIMATION OF THE ROTATION INDEX} 
An estimation of the rotation index of the closed curve $\alpha$ can be obtained by applying Cauchy-Crofton's formula \cite{SANTALO} in order to estimate the lengths of  $\alpha$ and the offset curve $\beta$, then applying Theorem \ref{MAIN}.

\vspace{2cm}

\noindent
{\small Enrique Mac\'{\i}as-Virg\'os\\
Institute of Mathematics\\
Department of Geometry and Topology\\
University of Santiago de Compostela\\
15782- SPAIN\\
{\tt xtquique@usc.es\\
\url{http://www.usc.es/imat/quique}\\
}

\end{document}